\documentclass[12pt]{article}
\usepackage{amsmath,amsfonts}
\usepackage{hyperref}
\usepackage{amsthm}
\usepackage{enumitem}
\usepackage{lettrine}
\usepackage{amssymb,amsmath,makeidx,verbatim}
\numberwithin{equation}{section}
\newtheorem{theorem}{Theorem}[section]
\newtheorem{corollary}[theorem]{Corollary}
\newtheorem{lemma}[theorem]{Lemma}
\newtheorem{remark}[theorem]{Remark}
\newtheorem{definition}[theorem]{Definition}
\newtheorem{proposition}[theorem]{Proposition}

\usepackage[british]{babel}

\usepackage{lettrine}
\date{}
\begin{document}
\title{Eigenvalues of Weighted-Laplacian under the extended Ricci flow}
\author{Abimbola Abolarinwa\thanks{Department of Mathematics and Statistics, Osun State College of Technology, P.M.B. 1011, Esa-Oke, Nigeria. Email:- A.Abolarinwa1@gmail.com}}
 
\maketitle
\begin{abstract}
Let $\Delta_\varphi = \Delta -\nabla \varphi \nabla$ be a symmetric diffusion operator with an invariant weighted volume measure $d\mu = e^{-\varphi} dv$ on an $n$-dimensional compact Riemannian manifold $(M,g)$, where $g=g(t)$  solves the extended Ricci flow.  In this article we study the evolution and monotonicty of the first nonzero eigenvalue of  $\Delta_\varphi$ and we obatin several monotone quantities along the  extended Ricci flow and its volume preserving version under some technical assumption. We also show that the eigenvalues diverge in a finite time for the case $n\geq 3$. Our results are natural extension of some known results for Laplace-Beltrami operator under various geometric flows. 
\\ \ \\
{\bf MSC (2010)}:  Primary 53C21, 53C44, Secondary 35P30, 58J35
\\ 
{\bf Keywords}: Witten-Laplacian, eigenvalues, Ricci flow, monotonicty,
curvature
\end{abstract}


\section{Introduction}

The present article discusses certain behaviours of the first nonzero eigenvalue $\Lambda =\Lambda(t)$ of the Witten-Laplacian on a compact Riemannian manifold $(M, g, d\mu)$, whose one-parameter family of metrics $g=g(t)$ solves the extended Ricci flow. We establish the evolution, monotonicity and finitely large-time behaviour of $\Lambda$ under the flow. We also find several monotone quantities under some technical assumption. In fact, we suppose the curvature of $M$ remains bounded throughout the flow existence time $t \in [0,T]$. Let $M$ be a connected, smooth, oriented, closed (compact without boundary) Manifold and $\varphi: M \to \mathbb{R}$ be a family of smooth functions. Then

\begin{definition}{\bf (Extended Ricci flow)}\\
The couple $(g(t),\varphi(t)), \ {t \in [0,T], T < \infty}$ is called the extended Ricci flow if it satisfies the following system of quasi-linear parabolic partial differential equations
 \begin{equation}\label{eq14}
\left. \begin{array}{l}
\displaystyle \frac{\partial}{\partial t} g(x,t) = - 2 Rc(x,t) + 2 \alpha(t) \nabla
\varphi(x,t) \otimes \nabla  \varphi(x,t)
\\ \ \\
\displaystyle \frac{\partial}{\partial t}  \varphi(x,t) = \Delta \varphi(x,t),
\end{array} \right.
\end{equation}
such that $(g(0),\varphi(0))= (g_0,\varphi_0)$. Here $Rc$ is the Ricci curvature of $M$,\ $\Delta$ is Laplace-Beltrami operator, $\nabla$ is the gradient operator and $\alpha(t)$ is a nonincreasing constant function of time, bounded below by $\alpha_n>0$ at all time.
\end{definition}
The coupled flow (\ref{eq14}) was first introduced by List in \cite{[Li08]}, where he studied the cases $\alpha_n=2$ and $\alpha_n=\frac{n-1}{n-2}$. He also proved its short-time existence on any compact Riemannian manifold. The motivation for this flow comes from General Relativity theory where the stationary solution of the flow corresponds to the solution of static Einstein Vacuum equation. Extended Ricci flow arises in a number of situations, e.g., when one studies Ricci flow on Warped product spaces \cite{[Tran],[Wil]}, when one couples Ricci flow with heat flow for harmonic maps \cite{[BaiTr],[Mu12]} and in quantum field theory (QFT) when one consider the so called renomalisation group (RG) flow \cite{[OSW]}.  In general (\ref{eq14}) behaves in a similar manner to Hamilton Ricci flow \cite{[Ha82]}.  In deed, if $\varphi$ is a constant then the flow (\ref{eq14}) degenerates to the Hamilton Ricci flow. Many researchers have extended the results on Ricci flow to the extended Ricci flow.
Owing to enormous roles of eigenvalues in understanding the geometry and topology of manifolds, the study of its behaviours  on time dependent metrics is a topic of current research, For examples see  \cite{[Ab15],[Ca07],[Ca08],[FXZ],[GPT],[Li07a],[Li07b],[Li10],[Lin1],[Ma06],[Zha12],[Zha13]}. An important application of eigenvalues occurs in the monotonicity of Perelman's $\mathcal{F}$-energy \cite{[Pe02]} in determining nonexistence of nontrivial steady or expanding breathers on compact manifolds which consequently led to the celebrated noncollapsing theorem and final resolution of the longstanding Poincar\'e conjecture. Following the standard notation in the study of the  extended Ricci flow we denote
$$\mathcal{S} := Rc - \alpha \nabla \otimes \nabla \varphi,$$
$$ \ {\cal S}_{ij} := R_{ij} - \alpha \varphi_i \varphi_j, $$ $$S=g^{ij}\mathcal{S}_{ij} := R_g-\alpha|\nabla \varphi|^2,$$
where $R_g$ is the scalar curvature of $(M,g)$. 

Let $M$ and $g(t)$ be as defined above we consider the first nonzero eigenvalue $\Lambda(t)$ of the Witten-Laplacian satisfying the eigenvalue problem
\begin{equation}
\Delta_\varphi f(x,t) = - \Lambda(t) f(x,t).
\end{equation} 
 Here
$$ \Delta_\varphi = \Delta_g-\nabla \varphi \cdot \nabla$$
is the so called Witten-Laplacian, which is a symmetric diffusion operator on $L^2(M,d\mu)$,
 $$\Delta_g = g^{ij}\Big( \frac{\partial^2}{\partial x^i \partial x^j} - \Gamma^k_{ij}\frac{\partial}{\partial x^k} \Big)$$ 
 is the Laplace-Beltrami operator and $\nabla$ is the Levi-Civita connection on $(M,g)$ with respect to the Christoffel's symbols $\Gamma^k_{ij}$. Let $d\nu$ be the Riemannian volume measure, the weighted volume measure $d\mu$ is given by
 $$d\mu = e^{\varphi(x)} d\nu$$
 where $\varphi \in C^\infty(M)$. In fact, the associated weighted measure $d\mu$ makes the Witten-Laplacian a self-adjoint operator.
 The following integration by parts formula
\begin{equation}
\int_M \Delta_\varphi u v d\mu = - \int_M \nabla u \nabla v d\mu  = \int_M u \Delta_\varphi v d\mu
\end{equation}
 holds for any $u,v \in C^\infty(M)$ with $M$ compact. When $\varphi$ is a constant function, the Witten-Laplacian is just the Laplace-Beltrami operator.
 
 \begin{definition}{\bf (Mini-Max principle)}\\
 The first non-zero eigenvalue $\Lambda_1(M,g,d\mu)$ is characterised as follows
\begin{equation}\label{eq33}
  \Lambda_1(t) = \inf_{0\neq f \in W_0^{1,2}(M) }\Bigg\{ \mathcal{E}(f,f)  \ \  : \ \int_M f^2d\mu =1, \ \int_Mfd\mu=0\Bigg\},
 \end{equation}
 where  $\mathcal{E}(f,u) = \int \langle \nabla f, \nabla u\rangle d\mu$ and $ W_0^{1, 2}(M,g,d\mu)$ is the completion of $C^\infty_0(M,g,d\mu)$ with respect to the norm
$$\|f\|_{W^{1,2}} = \Big(\int_M |f|^2 d\mu + \int_M | \nabla f|^2 d\mu \Big)^{\frac{1}{2}}.$$
\end{definition}

In what follows we assume that $\Lambda(t)$  exists and is $C^1$-differentiable under the extended Ricci flow deformation $(g(t),\varphi(t))$ in the given interval $0\leq t\leq T$.

The rest is planned as follows: Section \ref{sec2} contains statements of results and some relevant comments that will help to prove the results. In Section \ref{sec3} we highlight important facts and some backgroud results on the first eigenvalue of the Witten-Laplacian. Section \ref{sec4} discusses the proofs of Theorems \ref{thm13}, \ref{thm15} and \ref{thm17}, Corollary \ref{cor16} and Proposition \ref{prop25}. These results are extended to the case of normalised flow in Section \ref{sec5}. The last section discusses the behaviour of the spectrum as time grows along self-similar solution to the extended Ricci flow.. Here we found out that the spectrum diverges at a finite time under gradient Shrinking soliton.

\section{Statement of results}\label{sec2}
In this aricle we study the evolution of the first nonzero eigenvalue $\Lambda(t)$ of the Witten-Laplacian 
 and monotonicity of some quantities involving $\Lambda(t)$  under the extended Ricci flow. Let $f$ be corresponding normalised eigenfunction our main results are the following
\begin{theorem}\label{thm13}
Let $(M,g(t),\varphi(t),d\mu), \ t \in [0,T]$ be a solution to the extended flow (\ref{eq14}) on a closed manifold and $\Lambda(t)$ be the first nonzero eigenvalue of the Witten-Laplacian
  $\Delta_\varphi$ corresponding to the eigenfunction $ f(t, x)$, then $\Lambda(t)$ 
  evolves by 
\begin{equation}\label{eqn21}
 \frac{d}{d t} \Lambda(t) \geq  \Lambda(t) \int_M  	S f^2 \ d\mu 
-   \int_M S | \nabla f |^2  d\mu  + 2  \int_M \mathcal{S}_{ij} f_i f_j d\mu.  
\end{equation}	
\end{theorem}

The monotonicity of $\Lambda(t)$ here depends on the sign of  $S$. Note that just as the nonnegativity of the scalar curvature is preserved along the Ricci flow \cite{[Ha82]}, so also the nonnegativity of $S$ is preserved as long as the extended flow exists.

\begin{theorem}\label{thm15}
Let $(M,g(t),\varphi(t),d\mu), \ t \in [0,T]$ be a solution to the extended flow (\ref{eq14}) on a closed manifold. Let $\Lambda(t)$ be the first nonzero eigenvalue of the Witten-Laplacian
  $\Delta_\varphi$. Suppose
$$ S_(x,t)  \geq z(t)=  \frac{S_{min}(0)}{ 1 - \frac{2}{m} S_{min}(0)t }.$$
Then 
\begin{equation}
   \frac{d}{d t} \Big[ \Lambda(t) \exp \Big( - 2 \beta  \int_0^T z(t) dt \Big) \Big] \geq 0,
  \end{equation}
  where ${\cal S}_{ij} - \beta Sg_{ij} \geq 0$, \ $\beta \in (0,\frac{1}{n}]$ and 
  \begin{equation}
   \Lambda(t) \geq  \Lambda(0)  e^{ 2 \beta  \int_0^T   z(t) dt}
  \end{equation}
  for $0 \leq t \leq T$.
\end{theorem}

In some applications $S$ may be required to be a constant or  bounded by a constant.  In this situation we have the following  corollary.

\begin{corollary}\label{cor16}
Let $(M,g(t),\varphi(t),d\mu), \ t \in [0,T]$ be a solution to the extended flow (\ref{eq14}) on a closed manifold. Let $\Lambda(t)$ be the first nonzero eigenvalue of the Laplacian
 $\Delta_\varphi$. Then  if $S \geq C > 0$ in $M \times  [t_0, t ]$ for some uniform constant $C$ along the flow we have 
\begin{equation*}
   \frac{d}{d t}\log  \Lambda(t)  \geq  C 
\end{equation*}
    and 
\begin{equation*}
  \Lambda(t) \geq \Lambda(t_0) e^{C(t - t_0)}, \ \ \ \ \  \ t > t_0.
  \end{equation*}
\end{corollary}

Next we obtain a monotone quantity along the flow (\ref{eq14}).
\begin{theorem}\label{thm17}
With the assumption of Theorem \ref{thm15}.  The following quantity 
\begin{equation*}
\Lambda(t) \cdot \Big( z^{-1}_0 - \frac{2}{n}t \Big)^{\beta n}
\end{equation*}
is nondecreasing along the extended Ricci flow.
\end{theorem}

Recently it has been proved in \cite{[Li15]} that 
$$\limsup_{t\to T}\Big(\max_M R(t)\Big) = \infty$$
 and that $|\nabla \varphi |^2$ is uniformly bounded under the extended Ricci flow for the case $n \geq 3$ and $T<\infty$. Similarly, for the case $\varphi$ is a harmonic map between manifolds $:M \to N$, Muller \cite{[Mu12]} has shown that if $N$ has nonpositive sectional curvature the energy density of $\varphi$ is bounded by initial data. Observe that without additional assumption one can easily deduce that 
 $$\lim_{t\to T} S_{min}(t) = \infty.$$
Finally, we use the above estimate together with the extension of Reilly formula on compact manifold to prove that the eigenvalues of $\Delta_\varphi$ diverge as time $t$ increases.
\begin{proposition}\label{prop25}
Let $(M,g(t),\varphi(t),d\mu), \ t \in [0,T], T < \infty$ be a solution to the extended flow (\ref{eq14}) on a closed manifold. Let $\Lambda(t)$ be the first nonzero eigenvalue of 
 $\Delta_\varphi$. Then  
 \begin{equation}
 \lim_{t\to T} \Lambda(t) = \infty,
 \end{equation}
 where ${\cal S}_{ij} - \beta Sg_{ij} \geq 0,$ \ $\beta \in (0,\frac{1}{n}]$.
\end{proposition}
We equally demonstrate that $\lim_{t\to T} \Lambda(t) = \infty$ holds on gradient shrinking extended Ricci flow. See \ref{prop6} for details.

\section{The first eigenvalue of Witten-Laplacian}\label{sec3}
Let
$ \Delta_\varphi = \Delta_g-\nabla \varphi \cdot \nabla$
be the Witten-Laplacian, which is a symmetric diffusion operator on $L^2(M,g,d\mu)$, with weighted volume measure $d\mu = e^{\varphi(x)} d\nu$. A classical result tells us that $(M,g,d\mu)$ is a smooth metric  measure space and $\Delta_\varphi$ is self-adjoint and nonnegative definite with respect to the measure $d\mu$. Hence, the spectrum of the Witten-Laplacian is a set of an infinite sequence of points
\begin{align}\label{eq31}
0=\Lambda_0(\Delta_\varphi)<\Lambda_1(\Delta_\varphi) \leq \Lambda_2(\Delta_\varphi) \leq ...\Lambda_k(\Delta_\varphi) \leq ... \to \infty \ \ \ as\ \  k \to \infty
\end{align}
Thus, the set of all eigenvalues
of  $\Delta_\varphi$, counted with multiplicity, is an increasing sequence (\ref{eq31}).
The eigenvalue problem involving the Witten-Laplacian on a closed manifold consists in finding all possible real $\Lambda$ (eigenvalues) such that there exists non-trivial functions $f$ (eigenfunctions) satisfyng 
\begin{align}\label{eq32}
\Delta_\varphi f  = - \Lambda f.
\end{align} 
It can then be found that only constant functions correspond to $\Lambda_0=0$. The eigenfunctions are $L^2(M,g, d\mu)$ orthonormal basis $\{f_0, f_1, f_2,...,f_k,...\}$ of real $C^\infty(M,g,d\mu)$ function such that 
$$\Delta_\varphi f_j = - \Lambda_j f_j, \ \ \  j=1,2,...$$
while the eigenvalues are $L^2(M,g,d\mu)$ orthogonal. By the mini-max principle, the first non-zero eigenvalue $\Lambda_1(M,g,d\mu)$ can be characterised as follows
\begin{equation}\label{eq33}
  \Lambda_1(t) = \inf_{0\neq f \in W_0^{1,2}(M) }\Bigg\{ \mathcal{E}(f,f)  \ \  : \ \int_M f^2d\mu =1, \ \int_Mfd\mu=0\Bigg\},
 \end{equation}
 where  $\mathcal{E}(f,u) = \int \langle \nabla f, \nabla u\rangle d\mu$ and $ W_0^{1, 2}(M,g,d\mu)$ is the completion of $C^\infty_0(M,g,d\mu)$ with respect to the norm
$$\|f\|_{W^{1,2}} = \Big(\int_M |f|^2 d\mu + \int_M | \nabla f|^2 d\mu \Big)^{\frac{1}{2}}.$$
The first nonzero eigenvalue of $\Delta_\varphi$ is usually called the spectral gap in the literature and has been widely studied \cite{[BE],[CZ13],[Ma09]}. It satisfies
$$\Lambda_1(\Delta_\varphi) \geq C_0$$
for some positive constant $C_0$ under a lower boundedness condition on the Bakry-Emery Ricci tensor $Rc^\varphi$
$$Rc^\varphi = Rc_M + \nabla^2 \varphi \geq C_0g.$$
A key observation by Bakry and Emery in \cite{[BE]}  is the following weighted Bochner-Weitzenb\.ock formula
\begin{equation}\label{eq34}
  \frac{1}{2}\Delta_\varphi(|\nabla f|^2) = |\nabla^2f|^2+\langle\nabla f, \nabla(\Delta_\varphi f)\rangle + Rc^\varphi(\nabla f, \nabla f)
 \end{equation}
 An immediate consequence of formula (\ref{eq34}) is the lower bound estimate
 $$\Lambda \geq \frac{n(s+1)A}{(n(s+1)-1)}$$
 obtained by L. Ma \cite{[Ma09]} under an assumption on Bakry-Emery Ricci tensor
$$Rc^\varphi  \geq \frac{|\nabla \varphi|^2}{ns}+A$$
for some $A>0$ and $s>0$ (see \cite{[FLL13]} for related result). 

We remark that $Rc^\varphi$ is a natural extension of the Ricci tensor as derived in \cite{[BE]}. It arises naturally in the study of Ricci soliton \cite{[Pe02]}. We also remark that the discreteness of the spectrum of $\Delta_\varphi$ ensures the existence of solutions of eigenvalue problem (\ref{eq32}). The embedding of weighted Sobolev space $W^{1,2}_0 \hookrightarrow L^2$ is compact which is also equivalent to the discreteness of the spectrum of $(M,g,d\mu)$ \cite{[CZ13]} and \cite[Theorem 10.11]{[Gri09]}.

\section{The proofs of results}\label{sec4} 

   In this section, we consider the eigenvalues of the Witten-Laplacian under the extended Ricci flow (\ref{eq14}), assuming the least eigenvalue $\Lambda = \Lambda(t)$ is a function of time only. 
   \begin{lemma}\label{lem1}
Let $(M,g(t),\varphi(t),d\mu), \ t \in [0,T]$ be a solution to the extended flow (\ref{eq14}) on a closed manifold. Let $f\in C^\infty(M)$ be a smooth function on  $(M,g(t),d\mu)$ Then
\begin{eqnarray}\label{eq21}
&&\frac{\partial}{\partial t}S = \Delta_g S + 2|{\cal S}_{ij}|^2 + 2\alpha|\Delta_g \varphi|^2 \\
&&\frac{\partial}{\partial t} \Delta_\varphi f=  2 {\cal S}_{ij}f_{ij} - 2 {\cal S}_{ij}\varphi_i f_j + \Delta_\varphi f_t - (\varphi_t)_if_j - 2\alpha(\Delta \varphi)\varphi_i f_j.
\end{eqnarray}
\end{lemma}

\proof
The proofs follow from standard calculation, see for instance \cite{[Li08],[Mu12]} for the proof of (\ref{eq21}). Recall also from the above mentioned references that
$$\frac{\partial}{\partial t}g^{ij} = 2{\cal S}^{ij}, \  \frac{\partial}{\partial t} \Delta_g f = 2 {\cal S}_{ij}f_{ij} + \Delta_g f_t - 2\alpha(\Delta \varphi)\varphi_i f_j \ \ \mbox{and}\ \  \frac{\partial}{\partial t}d\mu = -S d\mu$$
Now for $\Delta_\varphi f = \Delta_g f - \nabla \varphi \cdot \nabla f$ we obtain 
$$\frac{\partial}{\partial t}\Big( \Delta_\varphi f\Big) = \frac{\partial}{\partial t}\Big(\Delta_g f \Big) - \frac{\partial}{\partial t}\Big(g^{ij}\varphi_if_j\Big)$$
where the result follows at once.

\qed

Next we discuss the proofs of Theorem \ref{thm13}.

\subsection*{Proof of Theorem \ref{thm13}} 
\proof  
   Let $M$ be a closed Riemannian manifold and $g(t)$ evolve by the extended Ricci flow (\ref{eq14}) in the interval $t \in [0,T]$. Let $f(t) = f(t, x)$ be the corresponding eigenfunction to the first nonzero eigenvalue $\Lambda(t) = \Lambda(t, f)$ of $ \Delta_\varphi$, i.e, 
  \begin{equation}\label{eq35}
  - \Delta_\varphi f(t, x) = \Lambda(t)f(t, x).
  \end{equation}
Taking derivative with respect to time, we have 
$$ - \Big(\frac{\partial}{\partial t} \Delta_\varphi \Big) f(t, x) - \Delta_\varphi   \frac{\partial}{\partial t} f (t, x) = \Big(\frac{d}{d t} \Lambda(t) \Big)f(t, x)  + \Lambda(t) \frac{\partial}{\partial t} f(t, x), $$
multiplying the above by $f(t, x)$ and integrate with respect to the weighted volume measure on $M$, we have 
$$ - \int_M f \Big(\frac{\partial}{\partial t} \Delta_\varphi \Big) f \ d\mu - \int_M f \Delta_\varphi   \frac{\partial}{\partial t}f \ d\mu  = \frac{d}{d t} \Lambda(t) \int_M   f^2 \ d\mu  + \Lambda(t) \int_M f  \frac{\partial}{\partial t} f \ d\mu.$$  
Notice that  by the application of integration by parts and (\ref{eq35})
$$ - \int_M \varphi \Delta   \frac{\partial}{\partial t} f \ d\mu  =  \lambda(t) \int_M f \frac{\partial}{\partial t} f d\mu,$$
then, by using the normalisation condition $\int_M f^2 d\mu=1$ we arrive at 
 \begin{equation}\label{eq36}
   \frac{d}{d t} \Lambda(t)  = -\int_M f \Big(\frac{\partial}{\partial t} \Delta_\varphi \Big) f \ d\mu.  
\end{equation}
Using the evolution of the Witten-Laplacian under the flow (\ref{eq14}) in Lemma \ref{lem1} we have 
\begin{align}\label{eq37}
  \frac{d}{d t} \Lambda(t)  & = - 2  \int_M  {\cal S}_{ij} f_{ij}f  d\mu + 2\alpha \int_M (\Delta \varphi)\varphi_i f_j f d\mu + 2 \int_M  {\cal S}_{ij} \varphi_i f_j f d\mu + \int_M (\varphi_t)_i f_j f d\mu.
  \end{align}  
We simplify the first and the last terms in the right hand side of (\ref{eq37}) further
\begin{align*}
- 2  \int_M \mathcal{S}_{ij}  f_{ij} f  d\mu & = - 2\int_M \mathcal{S}_{ij}  f_{ij} f e^{-\varphi} d\nu \\
\displaystyle & =  2 \int_M {\cal S}_{ij, i} f_j f d\mu 
+ 2 \int_M  \mathcal{S}_{ij}f_i f_j -2 \int_M {\cal S}_{ij}\varphi_if_j f  d\mu.
\end{align*} 
Note that by using the contracted second Bianchi identity $R_{ij,i} = 1/2R_i$ it is easy to obtain the following identty
\begin{equation}\label{eq211}
{\cal S}_{ij,i} = \frac{1}{2} S_{,i} -\alpha(\Delta\varphi)\varphi_i.
\end{equation} 
Using the identity(\ref{eq211}) we have  
\begin{align*}
2 \int_M {\cal S}_{ij, i} f_j f d\mu & =  \int _M  S_{,i} f_j f d\mu -2 \alpha \int_M(\Delta\varphi)\varphi_i f_j f d\mu \\
 \displaystyle & = - \int_M S (\Delta_\varphi f) f d\mu - \int_M S |\nabla f|^2 - 2 \alpha \int_M(\Delta\varphi f)\varphi_i f_j f d\mu 
 \end{align*}
 Similarly
 \begin{align*}
\int_M (\varphi_t)_i f_j f d\mu & = - \int_M \varphi_t ( f_j f e^{-\varphi})_j d\nu \\
\displaystyle & =  - \int_M \varphi_t(\Delta_\varphi f) f d\mu 
- \int_M \varphi_t |\nabla f|^2  d\mu \\
\displaystyle & =  \Lambda(t) \int_M \varphi_t f^2 d\mu 
- \int_M \varphi_t |\nabla f|^2  d\mu\\
\displaystyle & \geq 0.
\end{align*}
The last inequality follows from the claim of gradient estimate 
\begin{align}
|\nabla f|^2 \leq \Lambda(t)(f^2-1), \ \ \ \  \ \forall \  t \in [0,T].
\end{align}
which is very intrigue but interesting to prove and the fact that $\varphi$ satisfies the heat equation.

Putting all the above into (\ref{eq37}) we have 
$$   \frac{d}{d t} \Lambda(t) \geq  \Lambda(t) \int_M  	S f^2 \ d\mu 
-   \int_M S | \nabla f |^2  d\mu  + 2  \int_M \mathcal{S}_{ij} f_i f_j d\mu. $$
Hence, we have proved Theorem \ref{thm13} on the evolution of the first eigenvalue.

\qed

Before we prove Theorem \ref{thm15} we quickly make some comments about sign preservation on $S$ and about the condition $\mathcal{S}_{ij} - \beta S g_{ij} \geq 0$. To show that nonnegativity of $S$ is preserved along the flow, we use the evolution equation for $S(t)$ which is written as follows
$$\frac{d S}{dt} \geq \Delta S +\frac{2}{n} S^2$$
since $|{\cal S}_{ij}|^2 \geq \frac{1}{n}S^2$ by Cauchy-Schwarz inequality. Comparing $S$ with the solution of corresponding ODE
$$\frac{d z(t)}{dt} =\frac{2}{n} z^2(t), \ \ z(0) = z_0 = S_{min}(0)$$
by using the maximum principle argument at any time $t \in [0,T]$. Then
\begin{equation*}
S(t) \geq z(t) = \frac{1}{z^{-1}_0-\frac{2}{n}t}.
\end{equation*}
Clearly, $S_{min}(0)>0$ implies $S_{min}(t) \to \infty$ in finite time $T \leq n/2 S_{min}(0)<\infty$. This also implies that $R_{min}(t) \to \infty$ as $t \to T$ and thus $g(t)$ becomes singular in finite time $T_{sing}\leq T<\infty$.

By this we can prove the following Proposition using Hamilton maximum principle \cite[Theorem 9.1]{[Ha82]},\cite[Theorem 4.6]{[CK04]} for tensors.
\begin{proposition}
Let $g(t)$ be a smooth one parameter family of Riemannian metrics satisfying (\ref{eq14}). If 
$$(\mathcal{S}_{ij} - \beta S g_{ij})(x,0) \geq 0,$$
then
\begin{align}
(\mathcal{S}_{ij} - \beta S g_{ij})(x,t) \geq 0
\end{align}
for some $\beta \in [0,\frac{1}{n}]$ and all $t \in [0, T].$
\end{proposition}

\subsection*{Proof of Theorem \ref{thm15}}
\proof{}

By setting $ \mathcal{S}_{ij} - \beta S_g g_{ij} \geq 0$, $\beta \in (0,\frac{1}{n}]$, 
along the flow we have the following monotonicity formula from  Theorem \ref{thm13}
\begin{equation}\label{eq38}
 \frac{d}{d t} \Lambda(t) \geq \   \Lambda(t) \int_M S f^2 \ d\mu 
  + ( 2 \beta - 1) \int_M S | \nabla f |^2 \ d\mu.
\end{equation}
From the definition of $\Lambda(t)$ we have that 
\begin{equation}\label{eq39}
\Lambda(t) \int_M f^2 \ d\mu =  \int_M  | \nabla f |^2 \ d\mu
\end{equation}
and $\Lambda(t) > 0$. Suppose further that $S \geq S_{min}$ for all $t$, we know that 
\begin{equation}\label{eq310}
S(x,t) \geq z(t)
\end{equation}
Then using (\ref{eq39}) and (\ref{eq310}) we obtain  
\begin{equation}\label{eq311}
   \frac{d}{d t} \Lambda(t)  \geq  2 \beta z(t) \Lambda(t),
  \end{equation}
This concludes the proof of Thorem \ref{thm15}.

\qed

\subsection*{Proof of Corollary \ref{cor16}}
\proof
At any time $t \in[0,T]$, we write 
\begin{equation}
S(t) \geq \frac{1}{z^{-1}_0-\frac{2}{n}t} \geq \frac{n}{2|t+C_0|} = C
\end{equation}
for some constant $C_0$ and $C=C(g(0),t,n)$ depending on the geometry of the flow. The desired resultd follows at once.
\qed

\begin{remark}
If  $S$ is nonnegative then we are saying that the  curvature operator 
on the manifold is nonnegative. Precisely, $ R_g \geq \alpha |\nabla \varphi|^2$. 
Find more relevant result for this case in \cite{[Li10]}.
\end{remark}

\subsection*{Proof of Theorem \ref{thm17}}
\proof
Note that both $\Lambda_1(t)$ and $z(t)$ are functions of time only.
Denoting
$$Z(0) = S_{min}(0) = z_0,$$
we can evaluate
\begin{align*}
\int_{t_1}^{t_2} z(t) dt &= \int_{t_1}^{t_2} \Big(  \frac{1}{z^{-1}_0 - \frac{2}{n}t} \Big) dt = \log \Bigg( \frac{ z^{-1}_0 - \frac{2}{n}t_1}{ z^{-1}_0 - \frac{2}{n}t_2} \Bigg)^{\frac{n}{2}}.
\end{align*}
Therefore integrating both sides of (\ref{eq311}) from $t_1$ to $t_2$ yields 
\begin{equation}
\log  \frac{\Lambda_1(t_2)}{ \Lambda_1(t_2)} = \log \Bigg( \frac{ z^{-1}_0 - \frac{2}{n}t_1}{ z^{-1}_0 - \frac{2}{n}t_2} \Bigg)^{\beta n}
\end{equation}
for any time $t_1 < t_2$ sufficiently close to $t_2$. By this we have 
\begin{equation*}
\Lambda_1(t_2) \cdot \Big( z^{-1}_0 - \frac{2}{n}t_2 \Big)^{\beta n} \geq \Lambda_1(t_1) \cdot \Big(z^{-1}_0 - \frac{2}{n}t_1\Big)^{\beta n}.
\end{equation*}
Then $\Lambda_1(t) \cdot ( z^{-1}_0 - \frac{2}{n}t)^{\beta n}$ is nondecreasing along the extended Ricci flow. This completes the proof of Theorem \ref{thm17}.

\qed

Theorem \ref{thm17} has been proved for the eigenvalue of the $p$-Laplacian by the author in \cite{[Ab15],[AbMao]}.

Lastly in this section we want to apply the Reilly formula for Witten- Laplacian on a closed manifold
\begin{equation}\label{eq414}
\int_M \Big((\Delta_\varphi f)^2 - |\nabla^2f|^2\Big)d\mu = \int_M\Big(Rc+\nabla^2 \varphi\Big)\Big(\nabla f, \nabla f\Big)d\mu
\end{equation}
for $f\in C^\infty(M)$ to prove Proposition \ref{prop25}.

\subsection*{Proof of Proposition \ref{prop25}}
\proof
By the Reilly formaula extension (\ref{eq414}) we have 
\begin{equation}\label{eq415}
\int_M \Big((\Delta_\varphi f)^2 - |\nabla^2 f|^2\Big)d\mu = \int_M \mathcal{S}\Big(\nabla f, \nabla f\Big)d\mu +\int_M \Big(\alpha \nabla \varphi \otimes \nabla \varphi + \nabla^2\varphi\Big)\Big(\nabla f, \nabla f\Big)d\mu
\end{equation}
since $\Delta_\varphi f = \Delta f - \nabla \varphi \nabla f$.
We easily get the following inequality for $ s>0$
\begin{equation}\label{eq416}
(\Delta f)^2 = (\Delta_\varphi f + \nabla \varphi \nabla f)^2 \geq \frac{\Delta_\varphi f}{s+1} - \frac{|\nabla \varphi \nabla f|^2}{s}.
\end{equation}
By Cauchy-Schwarz inequality we have 
\begin{equation}\label{eq417}
|\nabla^2 f|^2 \geq \frac{1}{n}(\Delta f)^2 \geq \frac{\Delta_\varphi f}{n(s+1)} - \frac{|\nabla \varphi \nabla f|^2}{ns}
\end{equation}
Recall that 
$\Delta_\varphi f = -\Lambda f$
which implies
\begin{equation}\label{eq418}
\int_M\Big(\Delta_\varphi f \Big)^2 d\mu = \Lambda^2 \int_M f^2 d\mu. 
\end{equation}
Combinning (\ref{eq417}) and (\ref{eq418}) we obtain
\begin{equation}\label{eq419}
\int_M \Big((\Delta_\varphi f)^2 - |\nabla^2 f|^2\Big)d\mu \leq \Big(1-\frac{1}{n(s+1)}\Big) \Lambda^2 \int_M f^2 d\mu + \frac{1}{ns}\int_M \langle \varphi,\nabla f\rangle^2 d\mu
\end{equation}
Putting (\ref{eq419}) into (\ref{eq415}) yields
\begin{equation*}
\frac{n(s+1)-1}{n(s+1)}\Lambda^2 \int_M f^2 d\mu \geq  \int_M \mathcal{S}(\nabla f, \nabla f) d\mu + \int_M \nabla^2 \varphi (\nabla f, \nabla f) d\mu,
\end{equation*}
where we have used the facts that $\alpha \nabla \varphi \otimes \nabla \varphi (\nabla f \nabla f) \geq 0$  and that $|\nabla \varphi|^2$ is uniformly bounded for all $t \in [0,T]$. We now make use of the following observation 
$$|\nabla^2 \varphi| \geq \frac{1}{\sqrt{n}} |\Delta \varphi|^2 =  \frac{1}{\sqrt{n}}|\varphi_t|$$  
since $\varphi$ solves heat equation and the condition $\mathcal{S}_{ij} -\beta S g \geq 0$ to obtain
\begin{align*}
\displaystyle \frac{n(s+1)-1}{n(s+1)}\Lambda^2 \int_M f^2 d\mu &\geq \beta \int_M S|\nabla f|^2 d\mu + \frac{1}{\sqrt{n}} \int_M |\varphi_t| |\nabla f|^2 d\mu\\ \ \\
\displaystyle &\geq \Lambda \Big( \beta S_{min}(t)+  \frac{1}{\sqrt{n}} \min_M |\varphi_t |\Big)\int_M f^2 d\mu  
\end{align*}
and then
$$\Lambda(t) \geq \frac{n(s+1)}{n(s+1)-1}\Big( \beta S_{min}(t)+ \frac{1}{\sqrt{n}} \min_M |\varphi_t |\Big).$$
Taking the limit of both sides of the last inequality gives the result at once since 
 $ \lim_{t \to T} S_{min}(t) = \infty$
and $ \min |\varphi_t |$ is finite.

\qed
 
\section{Monotonicity under normalized flow}\label{sec5}
In this section we consider the volume preserving version of the extended Ricci flow. Let $(M,g(t),N,\varphi(t),d\mu)$ be as described before, the normalized extended Ricci flow is $(g(t),\varphi(t)), \ t \in [0,T], T\leq \infty$ satisfying
 \begin{equation}\label{eq41}
\left. \begin{array}{l}
\displaystyle \frac{\partial}{\partial t} g_{ij} = - 2\Big({\cal S}_{ij} -\frac{r}{n}g_{ij}\Big)
\\ \ \\
\displaystyle \frac{\partial}{\partial t}  \varphi = \Delta_g  \varphi,
\end{array} \right.
\end{equation}
on closed Riemannian manifold $M$, where $r=(Vol_{d\mu}(M))^{-1}(\int_M S d\mu)$ is the average of $S$. Notice that (\ref{eq14}) is a special case of (\ref{eq41}) for $r\equiv 0$.

 We only need to use the following evolution formulas instead of Lemma \ref{lem1}

\qed

\begin{lemma}\label{lem42}
Let $(M,g(t),\varphi(t),d\mu), \ t \in [0,T]$ be a solution to the extended flow (\ref{eq14}) on a closed manifold. Let $f\in C^\infty(M)$ be a smooth function on  $(M,g(t),d\mu)$ Then
\begin{align}\label{eq43}
\frac{\partial}{\partial t}S &= \Delta_g S + 2|{\cal S}_{ij}|^2 + 2\alpha|\Delta_g \varphi|^2 -\frac{2}{n}r S \\
\frac{\partial}{\partial t} \Delta_\varphi f &=  2 {\cal S}_{ij}f_{ij} - 2 {\cal S}_{ij}\varphi_i f_j + \Delta_\varphi f_t - (\varphi_t)_if_j - 2\alpha(\Delta \varphi)\varphi_i f_j -\frac{2}{n}r\Delta_\varphi f\\
\frac{d}{dt}d\mu &= (r-S) d\mu.
\end{align}
\end{lemma}

\begin{remark}
When $M$ is a compact Riemann surface, the Gauss-Bonnet Theorem asserts that 
$$r= \frac{4\pi \mathcal{X}}{Area(M^2)}$$
which is a constant, where $\mathcal{X}$ and $Area(M^2)$ are respectively the Euler characteristic  and area of the surface. In this case one can get many interesting monotone quantities depending on the sign of $r$. See for examples \cite[Corollary 4.2]{[FXZ]}, \cite[Corollary 1.5, Theorem 1.7]{[HL15]}, \cite[Section 3]{[CHL12]} and \cite{[Lin1]}.
\end{remark}

Finally, we state and prove our results for the eigenvalue $\Lambda_1(\Delta_\varphi)$ of the Witten-Laplacian under the normalized flow.

\begin{theorem}\label{thm44}
Let $(M,g(t),\varphi(t),d\mu), \ t \in [0,T], \ T\leq \infty$ be the solution to the normalized extended flow (\ref{eq41}) on a closed manifold. Let $\Lambda(t)$ be the first nonzero eigenvalue of the Witten-Laplacian
  $\Delta_\varphi$. Suppose
\begin{align}\label{eq46}
 S_(x,t) \geq y(t) = \frac{S_{min}(0)e^{-\frac{2}{n}\int_0^t r(s)ds}}{ 1 - \frac{2}{m} S_{min}(0)\int_0^t(e^{-\frac{2}{n}\int_0^t r(\tau)d\tau})ds }.
\end{align}
Then 
\begin{equation}
  \Lambda(t) e^{\frac{2}{n}\int_0^t r(s)ds}
  \end{equation}
  is nondecreasing along the flow (\ref{eq41}) and 
\begin{equation}
   \frac{d}{d t} \Big[ \Lambda(t) \exp \Big(2 \int_0^t (\frac{1}{n} r(s) - \beta y(s) )ds \Big) \Big] \geq 0,
  \end{equation}
or   equivalently
\begin{equation}
   \frac{d}{d t}\Big(\ln \Big( \Lambda(t) \exp \Big(\frac{2}{n} \int_0^t r(s) ds \Big)\Big) \geq 2 \beta y(t).
  \end{equation}
Moreover, we have the lower bound 
  \begin{equation}
   \Lambda(t)e^{\frac{2}{n} \int_0^t r(s) ds} \geq  \Lambda(0)  e^{ 2 \beta  \int_0^t y(s)ds}
  \end{equation}
  for \ $0 \leq t \leq T$.
\end{theorem}

\proof{}

Using the evolution of $\Delta_\varphi$ from Lemma \ref{lem42} into (\ref{eq36}), together with condition ${\cal S}_{ij} -\beta Sg_{ij}\geq 0$, \ $\beta \in (0,\frac{1}{n}]$ yields
\begin{equation}\label{eq412}
 \frac{d}{d t} \Lambda(t) \geq  -\frac{2}{n}\Lambda(t)+\Lambda(t) \int_M  	S f^2 \ d\mu 
 +(2\beta-1)   \int_M S | \nabla f |^2  d\mu.  
\end{equation}	
Recall the evolution of $S$ under normalised flow and the fact that $|\mathcal{S}_{ij}|^2 \geq \frac{1}{n}S^2$. Then
$$\frac{d S}{dt} \geq \Delta S +\frac{2}{n} S(S-r).$$
 Comparing $S$ with the solution of corresponding ODE
$$\frac{d y(t)}{dt} =\frac{2}{n} y(y-r), \ \ y(0) = y_0 = S_{min}(0)$$
by the maximum principle argument at any time $t \in [0,T]$ yields (\ref{eq46}). In particular when $S_{min}(0) =0$, then $S(0) \geq0$ and (\ref{eq46}) yields $S(x,t)\geq0$ for all $t\geq 0$. Then by (\ref{eq412}) together with the definiton of $\Lambda(t)$ we have
\begin{align*}
\frac{d \Lambda(t)}{dt} \geq \Big( 2\beta y(t) -\frac{2}{n} r \Big)\Lambda(t)
\end{align*}
and 
\begin{align*}
\frac{d}{dt} \Big(\ln \Lambda(t)\Big) \geq \Big( 2\beta y(t) -\frac{2}{n} r \Big)
\end{align*}
from where all the desired monotonicity and estimates follow by integrating from $0$ to $t$.

\qed

\section{Eigenvalues and extended Ricci soliton}\label{sec6}

Generally speaking,  a soliton is a  self-similar solution to an evolution equation which evolves along the symmetry group of the flow. In the case of the extended Ricci flow, the symmetries are scalings and diffeomorphisms. Soliton solutions are very crucial to the study of behaviour of solutions near singularities  of geometric flows where singularity models are more obvious.

\subsection*{Gradient Solitons} 
 In this case we modify the flow by a one-parameter group of diffeomorphisms $\psi_t$ and define a $time-$dependent 
 vector field $X_t$ from it.
\begin{definition}
 Let $\{ \psi_t \}, t \in I$ be a one-parameter  family of diffeomorphisms, 
 $\psi_t : M \rightarrow M,$ satisfying $\psi_0 = Id_M$ and  $\{g(t), \varphi(t)\}\ { t \in[0, T)}$ be a  one-parameter  family of solutions to the  extended Ricci flow 
 defined on $M$. Given a smooth scalar function $b(t) > 0 $, such that 
\begin{equation}\label{eq61}
g(t) = b(t) \psi^*_t g_0. \ \  {\mbox and} \ \ \varphi(t) = \psi_t^*\varphi(0).
 \end{equation} 
 Any pair  $\{g(t), \varphi(t)\}\ { t \in[0, T)}$ with this property is called an extended Ricci soliton.
 \end{definition} 
  This simply means that on an extended Ricci soliton all
  the Riemannian manifolds $(M^n, g)$ are isometric up to a scale factor that is allowed to vary 
  with $time$. Therefore, the extended Ricci soliton equation is equivalent to 
 \begin{equation}\label{eq62}
  \mathcal{S}_{ij}(g_0) + \frac{1}{2} \mathcal{L}_X g_0 = \sigma g_0 \ \ {\mbox and} \ \  (\Delta \varphi)( 0) = \mathcal{L}_X \varphi(0)
   \end{equation} 
    for any $ \sigma(t) = -\frac{1}{2} b'(t)$, where  $X$ is a vector field on $M$ and  $\mathcal{L}_X g_0$ is the Lie derivative of the evolving metric $g(t)$.
If the vector field $X$ is the gradient of a function, say $f$, then the solution is called a {\bf gradient  soliton}
and (\ref{eq62}) becomes a coupled elliptic system
 \begin{equation}\label{eq63}
 \left. \begin{array}{l}
\displaystyle   \mathcal{S}_{ij} + \nabla_i \nabla_j f =\sigma g_{ij}, \\ \ \\
 \displaystyle  \Delta \varphi -\langle \nabla \varphi, \nabla f \rangle = 0,  
 \end{array} \right. 
   \end{equation} 
where $\sigma$ is the homotheity constant. The case $b'(t) <0$, $b'(t)=0$ or $b'(t) >0$ corresponds to shrinking, steady or expanding gradient soliton. 
 \begin{lemma}
 Suppose the extended Ricci flow $(M, g(t), \varphi(t)), t \in [0, T]$, satisfies (\ref{eq61}), then, there exists a vector field $X$ on $M$ such that $(M, g_0, \varphi(0), X)$ 
 satisfies (\ref{eq63}). Conversely, given any solution $(M, g_0, \varphi(0), X)$ of  (\ref{eq63}), then, there exist a 
 one-parameter family of diffeomorphism $\psi_t$  of $M$  and scalar function $bt)$ 
 such that $(M, g(t), \varphi(t))$ of (\ref{eq61}) solves  the extended Ricci flow for all $ t \in [0, T]$.
 \end{lemma}

\begin{definition}
A $5$-turple $(M,g,\varphi,f,\sigma)$, where $f:M\to \mathrm{R}$ is a smooth function on $M$ and $\sigma \in \mathrm{R}$, is called a gradient extended Ricci soliton if it satisfies the coupled elliptic system
  \begin{equation}\label{eq64}
\left. \begin{array}{l}
\displaystyle \mathcal{S} + Hess f = \sigma g
\\ \ \\
\displaystyle \Delta  \varphi - \langle \nabla \varphi, \nabla f \rangle = 0,
\end{array} \right.
\end{equation}
where $ Hess f$ is the Hessian of potential function $f$. The soliton is said to be shrinking, steady or expanding if $\sigma>0, \ \sigma=0$ or $\sigma<0$ respectively.
\end{definition}
We say that the extended Ricci soliton is trivial if $f$ is a constant function, in which case $\varphi$ must be harmonic. Then (\ref{eq64}) reduces to 
 \begin{equation*}
\left. \begin{array}{l}
\displaystyle \mathcal{S} = \sigma g\ \ \ {\mbox and}\  \ \  \Delta  \varphi = 0
\end{array} \right.
\end{equation*}
which is a generalisation of Einstein manifold.
By M\"uller \cite{[Mu12]}, there exists a constant $C$ such that 
 \begin{equation}\label{eq65}
 S+|\nabla f|^2 - 2 \sigma f = C.
\end{equation}
On the other hand, taking the trace of the first equation in (\ref{eq64}) we have 
\begin{equation}\label{eq66}
 S+ \Delta f = n \sigma.
\end{equation}
By combinning (\ref{eq65}) and (\ref{eq66}) we obtain
\begin{equation}\label{eq67}
 -\Delta f +|\nabla f|^2 - 2 \sigma f =(C- n \sigma)
\end{equation}
which implies
\begin{equation}\label{eq68}
 -\Delta_f f - 2 \sigma f =C_0,
\end{equation}
 where $C_0=C- n \sigma$.
 
A further analysis shows that $2\sigma$ is an eigenvalue of $\Delta_f$ by allowing $f$ to satisfy $\int_M fe^{-f} d\mu =0$.
\begin{lemma}
$2\sigma$ is an eigenvalue of $\Delta_f$, that is, 
$$ \Delta_f f = - 2\sigma f.$$
\end{lemma}
 
The solutions defined by (\ref{eq64}) are special solutions for the extended Ricci flow \cite{[Li08],[Mu12]}. Notice that if $\varphi : M \to \mathrm{R}$ is a constant then  (\ref{eq64}) reduces to 
\begin{equation}\label{eq69}
 Rc + Hess f = \sigma g
\end{equation} 
the gradient Ricci soliton, a special solution of Ricci flow. Since the extended Ricci flow is a generalisation, it is natural to ask whether the same results under the Ricci flow hold under extended Ricci flow.

\begin{theorem} (\cite[Theorem 3.1]{[GPT15]})
\\
Suppose $M$ is compact and $(M,g,\varphi,f,\sigma)$, is a gradient extended Ricci soliton. Then $S_{min}=0$ in the steady case,  $0\leq S_{min}\leq -n\sigma$ in the shrinking case and in particular $S\geq 0$; in the expanding case $-n\sigma \leq S_{min} \leq 0$ and in particular $S \geq -n \sigma$.
\end{theorem}
A corolary of the above is that $\mathcal{S}_{ij} = -\sigma g$ on a steady or an expanding soliton given that potential function $f$ is constant or sub-harmonic and $\varphi$ is a harmonic map (see \cite[Corollary 3.2]{[GPT15]}).

In conclusion of this section we determine the behaviour of the spectrum of $M$ on special solutions; namely the gradient shrinking and expanding extended Ricci flow.
 
\begin{proposition}\label{prop6}
Let $(M,g(t),\varphi(t),d\mu), \ t \in [0,T], T < \infty$ be a gradient shrinking soliton. Let $\Lambda(t)$ be the first nonzero eigenvalue of 
 $\Delta_\varphi$ with $f$ being the associated normalised eigenfunction. Then  
 \begin{equation}
 \lim_{t\to T} \Lambda(t) = \infty
 \end{equation}
 where ${\cal S}_{ij} - \beta Sg_{ij} \geq 0$, \ $\beta \in (0,\frac{1}{n}]$.
\end{proposition} 

\proof
By a straightforward calculation  using (\ref{eq66}) we have
 \begin{equation}\label{eq611}
\left. \begin{array}{l}
\displaystyle \int_M (\Delta f)^2 d\mu= -\int_M S \Delta f d\mu 
= \int_M \nabla f \nabla(S e^{-\varphi}) dv\\
\displaystyle \hspace{2.2cm} = \int_M \langle \nabla f, \nabla S\rangle d\mu - \int_M \langle \nabla \varphi, \nabla f\rangle S d\mu \\
\displaystyle \hspace{2.2cm} = 2 \int_M \mathcal{S}(\nabla f, \nabla f)d\mu- \int_M \langle \nabla \varphi, \nabla f\rangle S d\mu
\end{array} \right.
\end{equation}
where we have used the following identity 
$$\nabla S = 2 \mathcal{S}(\nabla f, \cdot).$$
This identity follows from a standard computation in the theory of extended Ricci soliton (See \cite[Lemma 2.1]{[GPT15]}).
Observe that
\begin{align}\label{eq612}
(\Delta f)^2 = (\Delta_\varphi f)^2 + \langle \nabla \varphi, \nabla f\rangle^2 + 2 \Delta_\varphi \langle \nabla \varphi, \nabla f\rangle
\end{align}
and that the following inequality
\begin{align}\label{eq613}
 - \Delta_\varphi \langle \nabla \varphi, \nabla f\rangle \geq - \frac{(\Delta_\varphi f)^2}{2a} +\frac{a \ \langle \nabla \varphi, \nabla f\rangle^2}{2}
 \end{align}
 holds for $a \geq 1.$
 Recall also that in the case of gradient shrinking soliton $-S\geq  n\sigma$. 
 
  We now put (\ref{eq612}) and (\ref{eq613}) into (\ref{eq611}) and  obtain
\begin{align}\label{eq614}
\left. \begin{array}{l}
\displaystyle
\Big(1 +\frac{1}{a}\Big)\int_M (\Delta_\varphi f)^2 d\mu \geq 2 \int_M \mathcal{S}(\nabla f, \nabla f)d\mu +n\sigma \int_M \langle \nabla \varphi, \nabla f\rangle d\mu \\ \ \\
\displaystyle \hspace{5cm} +(a-1)\int_M \langle \nabla \varphi, \nabla f\rangle^2 d\mu\\ \ \\
\displaystyle \hspace{4.2cm} \geq 2 \beta \int_M S|\nabla f|^2 d\mu +n\sigma \int_M \Delta \varphi d\mu
\end{array} \right.
\end{align} 
since gradient soliton implies $\Delta \varphi - \langle \nabla \varphi, \nabla f\rangle =0$ and $a \geq 1$. 

Note that $\int_M \Delta \varphi d\mu =0$ on a closed manifold, then we have 
\begin{equation}
\Lambda(t) \geq \frac{2a \beta}{a+1} \max_M S_{min}(t)
\end{equation}
Taking the limit as $t$ grows to maximal time $T$ gives the desired result.

\qed
\begin{remark}
In the above case the quantity $\int_M \mathcal{S}(\nabla f, \nabla f)d\mu$ is positive. If it was non-positive then the soliton must be trivial. In fact, compact shrinking soliton has been characterised to be trivial in \cite[Theorem 3.2]{[Tad]} if and only if 
\begin{align*}
\int_M \mathcal{S}(\nabla f, \nabla f)d\nu \leq \frac{\lambda_1}{2} \int_M \|\nabla f\|^2 d\nu
\end{align*}
where $\lambda_1$ is the first non-zero eigenvalue of the Laplacian.
\end{remark}

Each soliton can be written in a canonical form by choosing $b(t)$ in (\ref{eq61}) to be 
$$b(t) = 1+2\sigma t.$$
We also have that a diffeomorphism $\psi: (M, \psi^*g,d\mu) \to (M, g,d\mu)$ is an isometry, which implies that $(M, \psi^*g,d\mu)$ and $(M, g,d\mu)$ have the same spectrum with eigenfunction $\{f_j\}$ and $\{\psi^*f_j\}$, $j = 1,2,...$, respectively. Now if $(M,g(t),\varphi(t), d\mu)$ is an extended Ricci soliton on $(M, g(0), \varphi(0),d\mu)$ then, the spectrum $spc(g(t))$ at time $t$ is proportional to the initial spectrum $spc(g(0))$, that is,
\begin{equation}
spc(g(t)) = \frac{1}{1+2\sigma t} spc(g(0)).
\end{equation}
We then  have 
\begin{align}\label{eq617}
\displaystyle \frac{d\Lambda(t)}{dt} = - \frac{2\sigma}{(1+2\sigma t)^2}
\end{align}
as the change rate of the spectrum of $(M, g(t), d\mu)$. If we combine (\ref{eq617}) and (\ref{eqn21}) of Theorem \ref{thm13} with the conditon ${\cal S}_{ij} - \beta S g \geq 0$, \ $\beta \in (0,\frac{1}{n}]$,  we get the following 
\begin{align*}
\displaystyle - \frac{2\sigma}{(1+2\sigma t)^2} \geq \Lambda(t) \int_M S f^2 \ d\mu 
  + ( 2 \beta - 1) \int_M S | \nabla f |^2 \ d\mu.
\end{align*}
Since we have $S \geq -n\sigma$ in the expanding case, then
\begin{align*}
\displaystyle - \frac{2\sigma}{(1+2\sigma t)^2} \geq -2\beta n \sigma \Lambda(t).
\end{align*}
This yields the following lower bound for the first eigenvalue 
\begin{align*}
\displaystyle  \Lambda(t) \geq  \frac{1}{n\beta(1+2\sigma t)^2}.
\end{align*}
on the expanding soliton.
\begin{proposition}\label{prop7}
Let $(M,g(t),\varphi(t),d\mu), \ t \in [0,T], T < \infty$ be a gradient expanding soliton. Let $\Lambda(t)$ be the first nonzero eigenvalue of 
 $\Delta_\varphi$ with $f$ being the associated normalised eigenfunction. Then  
 \begin{align*}
\displaystyle  \Lambda(t) \geq  \frac{1}{n\beta(1+2\sigma t)^2}.
\end{align*}
 where ${\cal S}_{ij} - \beta Sg_{ij} \geq 0$ \ $\beta \in (0,\frac{1}{n}]$.
\end{proposition}

\subsection*{Acknowledgement} 
This research was carried out when the author was in AIMS-Cameroon on tutor fellowship between September--December, 2015. He is therefore grateful to  the director of the centre, Prof. Mama Foupouagnigni, for the hospitality during this period.


\end{document}